\documentclass[psamsfonts, reqno, 12pt]{amsart}

\usepackage{amsmath,amsfonts,amsthm,amssymb,psfrag,graphicx}

\DeclareMathOperator\ord{ord}

\numberwithin{equation}{section}

\newtheorem{theorem}{Theorem}[section]

\newtheorem{proposition}[theorem]{Proposition}

\theoremstyle{definition}
\newtheorem{example}[theorem]{Example}

\newtheorem{remark}[theorem]{Remark}

\newcommand{\fix}{F}

\newcommand{\ent}{{h}}
\newcommand{\mah}{\mathsf{m}}

\newcommand{\ch}{\mathrm{char}}
\newcommand{\ass}{\mathrm{Asc}}

\title[Uniform Periodic Point Growth in Entropy Rank One]{Uniform Periodic Point Growth\\in Entropy Rank One}
\author{Richard Miles}
\author{Thomas Ward}
\date{\today}

\address{School of Mathematics, University of East Anglia,
Norwich, NR4 7TJ, UK}
\subjclass[2000]{22D40, 37A15, 37A35}
\thanks{This research was supported by E.P.S.R.C. grant EP/C015754/1;
both authors express their thanks to Graham Everest and
Shaun Stevens for helpful discussions.}

\begin{document}

\begin{abstract}
We show that algebraic dynamical systems with entropy rank one have
uniformly exponentially many periodic points in all directions.
\end{abstract}
\maketitle

\section{Introduction}

Let~$\alpha$ be an action of~$\mathbb Z^d$ by
continuous automorphisms of a
compact metrizable abelian group~$X$ (such a
system is called an algebraic~$\mathbb Z^d$-action). For a continuous map~$\beta:X\to X$
write~$\ent(\beta)$ for the topological entropy and~$\fix(\beta)=\{x\in X\mid\beta x=x\}$
for the set of fixed points.
The action~$\alpha$ is said to have
\emph{entropy rank one} if, for
each~$\mathbf n\in\mathbb Z^d$,~$\ent(\alpha^{\mathbf n})<\infty$.

If~$\alpha$ is a mixing entropy rank one action and the topological
dimension~$\dim(X)$ is finite,~$\fix(\alpha^{\mathbf n})$ is finite
for all~$\mathbf n\in\mathbb Z^d$. Our purpose here is to show that
under natural conditions~$\vert\fix(\alpha^{\mathbf n})\vert$
exhibits uniform exponential growth. Write~$\mathbf n_j\to\infty$
for a sequence in~$\mathbb Z^d$ if for any finite
set~$F\subset\mathbb Z^d$ there is some~$J$ for which~$j>J$ implies
that~$\mathbf n_j\notin F$. Equivalently, this means~$\Vert\mathbf
n_j\Vert\to\infty$ as~$ j\to\infty$ where~$\Vert\cdot\Vert$ is the
Euclidean norm on~$\mathbb Z^d$. The Noetherian condition mentioned
in Theorem~\ref{theorem:uniformexponentialgrowth} is explained in
Section~\ref{section:mainsection}.

\begin{theorem}\label{theorem:uniformexponentialgrowth}
Let~$\alpha$ be a mixing algebraic~$\mathbb Z^d$-action with entropy
rank one on a finite-dimensional group $X$. Then there exist
constants $C_1,C_2\geq 0$ such that
\[
\limsup_{\mathbf n\to\infty}\frac{1}{\Vert\mathbf n
\Vert}\log\vert\fix(\alpha^{\mathbf n})\vert=C_1<\infty
\]
and
\[
\liminf_{\mathbf n\to\infty}\frac{1}{\Vert\mathbf n
\Vert}\log\vert\fix(\alpha^{\mathbf n})\vert=C_2.
\]
If~$\dim(X)>0$ and the action is Noetherian then~$C_2>0$.
\end{theorem}

An affirmative answer to Lehmer's problem
(see Remark \ref{remark:lehmersproblem}) would render
the finite dimension assumption in
Theorem~\ref{theorem:uniformexponentialgrowth} redundant.
The assumption is also not required for a mixing
Noetherian entropy rank one algebraic action; in
that setting it is a consequence of the Noetherian condition.
Essentially, the only
non-trivial conclusion in the theorem is that~$C_2>0$.

The action~$\alpha$ is \emph{expansive} if there
exists a neighbourhood $U$ of $0_X$ such that
\[
\bigcap_{\mathbf{n}\in \mathbb{Z}^d}\alpha^{\mathbf{n}}(U) = \{0_X\}.
\]
For~$d=1$ and~$\alpha$ expansive, the growth rate of periodic points
exists, so~$C_1=C_2$ in
Theorem~\ref{theorem:uniformexponentialgrowth}. The constant
coincides with the entropy and is non-zero provided~$X$ is infinite.
For~$d>1$, the zero-dimensional
Example~\ref{example:ledrappiersexample} is expansive yet
has~$C_2=0$. Expansive actions on connected groups are more
indicative of the import of a positive value for~$C_2$. In terms of
expansive subdynamics (see~\cite{MR1355295} and~\cite{MR1869066}),
there are sequences~$\mathbf n\to\infty$ converging to non-expansive
lines; along such sequences~$\vert\fix (\alpha^{\mathbf n})\vert$ is
much smaller than the same expression with~$\mathbf n$ of similar
Euclidean size and far from non-expansive directions. Nonetheless,
there is a uniform exponential growth in all directions.

\begin{example}
\label{example:firstvisittox2x3example} Consider
the~$\mathbb{Z}^2$-action~$\alpha$ dual to the~$\mathbb Z^2$-action
generated by the commuting maps~$r\mapsto 2r$ and~$r\mapsto 3r$
on~$\mathbb{Z}[\frac{1}{6}]$. Figure~\ref{figurex2x3periodicpoints}
shows the map~$\mathbf n\mapsto\vert\fix(\alpha^{\mathbf n})
\vert\in\mathbb N$. Notice that~$\vert\fix(\alpha^{-\mathbf
n})\vert=\vert\fix(\alpha^{\mathbf n})\vert$, so only the
region~$n_2\ge0$ is shown, with~$\infty$ denoting the lattice
point~$(0,0)$ corresponding to the identity map~$\alpha^{(0,0)}$.

\begin{figure}[ht]
\[
\begin{array}{ccccccccccc}
{211}&227&235&239&241&121&485&971&1943&3887&{\it 7775}\\
49&65&73&77&79&5&161&323&647&{\it 1295}&2591\\
{\bf 5}&11&19&23&25&13&53&107&{\it 215}&431&863\\
23&{\bf 7}&{\bf 1}&5&7&1&17&{\it 35}&71&143&287\\
29&13&{\bf 5}&1&1&1&{\it 5}&11&23&47&95\\
31&5&7&{\bf 1}&{\bf 1}&\infty&1&1&7&5&31\\
\end{array}
\]
\caption{\label{figurex2x3periodicpoints}Periodic
point counts
for~$\times2,\times3$.}
\end{figure}
In an expansive direction like~$(1,1)$ we
have~$\vert\fix(\alpha^{n(1,1)})\vert=6^n-1$ and the exponential
growth rate along the sequence in italics is clear. For lattice
points close to the non-expansive line~$2^x3^y=1$ we find (for
example) that~$\vert \fix(\alpha^{(-5,3)})\vert=5$, and it is not
immediately clear that there is exponential growth along the
sequence shown in bold;
Theorem~\ref{theorem:uniformexponentialgrowth} asserts that there
is.
\end{example}

\begin{example}\label{example:ledrappiersexample}
Consider the~$\mathbb{Z}^2$-action~$\alpha$ on~$X$ dual to
the~$\mathbb Z^2$-action generated by the commuting maps~$r\mapsto
u_1r$ and~$r\mapsto u_2r$ on
\[
\mathbb{Z}[u_1^{\pm1},u_2^{\pm1}]/\langle 2,1+u_1+u_2\rangle
\]
(this is an example of the type introduced by
Ledrappier~\cite{MR512106};~$X$ is a zero-dimensional group). As
shown
in~\cite{MR1461206},~$\vert\fix(\alpha^{(n,0)})\vert=2^{n-2^{\ord_2(n)}}$.
In particular,
\[
\lim_{n\to\infty}\frac{1}{\Vert(2^n,0)\Vert}\log\vert\fix(\alpha^{(2^n,0)})\vert=0,
\]
showing that some assumption on the topological dimension on~$X$ is
needed to have~$C_2>0$ in
Theorem~\ref{theorem:uniformexponentialgrowth}.
\end{example}

\begin{example}
The Noetherian condition is needed to have~$C_2>0$ even for~$d=1$.
For example, the automorphism~$\alpha$ dual to the map~$r\mapsto 2r$
on~$\mathbb Q$ has~$\vert\fix(\alpha^n)\vert=1$ for all~$n$. For
rings between~$\mathbb Z$ and~$\mathbb Q$ a variety of exotic
periodic point behavior is possible (see~\cite{MR1458718}
or~\cite{MR1619569} for the details).
\end{example}

Theorem~\ref{theorem:uniformexponentialgrowth} does however apply to
Noetherian non-expansive systems. Thus, for example, the genuinely
partially hyperbolic systems like that described by Damjanovi{\'c}
and Katok~\cite[Ex.~7.3]{dkpreprint} satisfy the hypotheses.

\section{Proof of Theorem~\ref{theorem:uniformexponentialgrowth}}\label{section:mainsection}

The case~$d=1$ is covered by~\cite{MR1461206} and~\cite{miles_zeta},
so we may assume~$d\ge2$. Algebraic~$\mathbb Z^d$-actions have a
convenient description in terms of commutative algebra due to
Kitchens and Schmidt~\cite{MR1036904}.

Let~$R_d=\mathbb{Z}[u_1^{\pm 1}, \dots, u_d^{\pm 1}]$ be the ring of
Laurent polynomials in commuting variables~$u_1, \dots, u_d$ with
integer coefficients. If~$\alpha$ is an algebraic~$\mathbb
Z^d$-action of a compact abelian group~$X$, the character
group~$M=\widehat{X}$ has the structure of a discrete
countable~$R_d$-module, obtained by first identifying the dual
automorphism~$\widehat{\alpha}^{\mathbf{n}}$ with multiplication by
the monomial~$u^{\mathbf{n}}=u_1^{n_1}\dots u_d^{n_d}$, and then
extending additively to multiplication by polynomials. Conversely,
any countable~$R_d$-module~$M$ has an associated
algebraic~$\mathbb{Z}^d$-action obtained by dualizing the action
induced by multiplying by monomials on~$M$. The
action~$\alpha=\alpha_M$ is described as \emph{Noetherian} if the
module~$M$ is Noetherian. A full account of this correspondence and
the resulting theory is given in Schmidt's
monograph~\cite{MR1345152}.

Entropy rank one actions are described in \cite{MR2031042} and
developments concerning their periodic points may be found in the
papers~\cite{miles_zeta} and~\cite{mileswarddetects}. By a slight
abuse of notation, write~$\ent(\cdot)$ for the topological entropy
of maps and for the extension of the entropy function to all
of~$\mathbb R^d$ in the sense explained below.
\begin{proposition}
Let~$\alpha_M$ be a mixing algebraic~$\mathbb Z^d$-action
with entropy rank one on a finite-dimensional group. Then
\begin{enumerate}
\item The set of associated primes of the associated~$R_d$-module,~$\ass(M)$ is
finite. For each~$\mathfrak p\in\ass(M)$, the domain~$R_d/\mathfrak
p$ has Krull dimension~$1$ and its field of fractions~$\mathbb
K(\mathfrak p)$ is a global field.
\item There exist positive integers~$m(\mathfrak p)$,~$\mathfrak p\in\ass(M)$, such
that for every non-zero~$\mathbf n\in\mathbb Z^d$,
\begin{equation}\label{equation:generalmoduleentropy}
\ent(\alpha_M^{\mathbf n})
=\sum_{\mathfrak p\in\ass(M)}m(\mathfrak p)\ent(\alpha_{R_d/\mathfrak p}^{\mathbf n})
\end{equation}
and
\begin{equation}\label{equation:generalmoduleperiodicpoints}
\vert\fix (\alpha_M^{\mathbf n})\vert
\leqslant
\prod_{\mathfrak p\in\ass(M)}
\vert\fix (\alpha_{R_d/\mathfrak p}^{\mathbf n})\vert^{m(\mathfrak p)},
\end{equation}
with equality if~$M$ is Noetherian.
\item For each~$\mathfrak p\in\ass(M)$, there is a finite set of
places~$\mathcal S(\mathfrak p)$ of~$\mathbb K(\mathfrak p)$ such
that
\begin{equation}\label{equation:primemoduleentropy}
\ent(\alpha_{R_d/\mathfrak p}^{\mathbf n})
=
\sum_{v\in\mathcal S(\mathfrak p)} \max\{\mathbf l_v\cdot\mathbf n,0\} > 0,
\end{equation}
and
\begin{equation}\label{equation:primemoduleperiodicpoints}
\vert\fix (\alpha_{R_d/\mathfrak p}^{\mathbf n})\vert
=\prod_{v\in\mathcal S(\mathfrak p)}\vert{\xi}_{\mathfrak
p}^{\mathbf n} - 1\vert_v,
\end{equation}
where~${\xi}_{\mathfrak p}=(\overline{u}_1,\dots,\overline{u}_d)$
denotes the image of~$u$ in~$(R_d/\mathfrak p)^d$ and~$\mathbf
l_v=(\log|\overline{u}_1|_v,\dots,\log|\overline{u}_d|_v)$.
\end{enumerate}
\end{proposition}

\begin{proof}
Let~$\mathbf n\in\mathbb Z^d$ be non-zero. Since~$\alpha_M$ is
mixing,~\cite[Prop. 6.6]{MR1345152} shows that for each~$\mathfrak
p\in\ass(M)$,~$\alpha_{R_d/\mathfrak p}^{\mathbf n}$ is ergodic,
so~$\ent(\alpha_{R_d/\mathfrak p}^{\mathbf n})>0$. It follows
from~\cite{MR2031042} that for each~$\mathfrak
p\in\ass(M)$,~$R_d/\mathfrak p$ has Krull dimension~$1$ and~$\mathbb
K(\mathfrak p)$ is a global field. If~$\ch(R_d/\mathfrak p)>0$
then~$\ent(\alpha_{R_d/\mathfrak p}^{\mathbf n})\ge\log 2$. Via
Yuzvinski{\u\i}'s formula (see~\cite[Th. 14.1]{MR1345152}
or~\cite{MR0214726}), this implies there are only finitely many
such~$\mathfrak p\in\ass(M)$, as~$\ent(\alpha_{R_d/\mathfrak
p}^{\mathbf n})<\infty$. Also, there can be only finitely
many~$\mathfrak p\in\ass(M)$ with~$\ch(R_d/\mathfrak p)=0$, as
\[
\dim(X)\geqslant\vert\{\mathfrak p\in\ass(M):\ch(R_d/\mathfrak p)=0\}\vert.
\]
This establishes~(1).

The method of~\cite[Lem.~8.2]{MR2031042} shows that in any prime
filtration of a Noetherian submodule of~$M$, each prime~$\mathfrak
p\in\ass(M)$ appears with a maximum multiplicity
\[
m(\mathfrak p)=\dim_{\mathbb K(\mathfrak p)}(M\otimes_{R_d}\mathbb K(\mathfrak p)),
\]
which is finite by similar reasoning to the proof of~(1). By
adopting an algorithm for obtaining a prime filtration which selects
the associated primes of~$M$ first, one obtains a Noetherian
submodule~$L\subset M$ such that each prime~$\mathfrak p\in\ass(M)$
appears with multiplicity~$m(\mathfrak p)$ in a filtration of~$L$.
Furthermore, if~$L\neq M$, each~$\mathfrak q\in\ass(M/L)$ is maximal
and~$R_d/\mathfrak q$ is a finite field. Hence, Yuzvinski{\u\i}'s
formula shows that~$\ent(\alpha_{M/L}^{\mathbf n})=0$
and~$\ent(\alpha_{M}^{\mathbf n})=\ent(\alpha_{L}^{\mathbf n})$; the
formula~\eqref{equation:generalmoduleentropy} then follows
from~\cite[Lem.~4.3]{MR1062797}.

If~$M$ is Noetherian, equality
in~\eqref{equation:generalmoduleperiodicpoints} is given
by~\cite[Th.~3.2]{miles_zeta}. If~$M$ is not
Noetherian,~\eqref{equation:generalmoduleperiodicpoints} follows
from~\cite[Th.~3.2]{miles_zeta} applied to~$L$, together with the
inequality
\[
\vert\fix (\alpha_M^{\mathbf n})\vert
\leqslant
\vert\fix (\alpha_L^{\mathbf n})\vert,
\]
which is established using a similar method to the proof
of~\cite[Lem.~2.6]{miles_periodic_points}.

Finally, the entropy formula~\eqref{equation:primemoduleentropy}
follows from~\cite[Prop.~8.5]{MR2031042} and the periodic point
counting formula~\eqref{equation:primemoduleperiodicpoints}
is~\cite[Lem~3.1]{miles_zeta}.
\end{proof}

\begin{proof}[Proof of Theorem \ref{theorem:uniformexponentialgrowth}]
Let~$M=\widehat{X}$ denote the dual~$R_d$-module and let~$\mathfrak
p\in\ass(M)$. For a fixed~$\mathfrak p\in\ass(M)$ and any
non-zero~$\mathbf n\in\mathbb Z^d$, set
\begin{equation}\label{equation:normalizedprimeperiodicpoints}
f(\mathbf n)=\frac{1}{\Vert\mathbf n
\Vert}\log\vert\fix(\alpha_{R_d/\mathfrak p}^{\mathbf n})\vert.
\end{equation}
Let~$h:\mathbb R^d\to\mathbb R_{\geqslant 0}$ denote the directional
entropy function for the action~$\alpha_{R_d/\mathfrak p}$. This is
the function obtained by extending the entropy
formula~\eqref{equation:primemoduleentropy} to values of~$\mathbb
R^d$ (see~\cite[Sec. 8]{MR2031042} for further details).

For any vector~$\mathbf v\in\mathbb R^d\setminus\{0\}$,
write~$\widehat{\mathbf v}\in\mathsf S_{d-1}$ for the unique vector
of unit length with the property that~$\mathbf
v=\lambda\widehat{\mathbf v}$ for some scalar~$\lambda>0$.
From~\eqref{equation:primemoduleentropy}
and~\eqref{equation:primemoduleperiodicpoints} it follows that
\begin{equation}\label{equation:decomposef}
f(\mathbf n)=g(\mathbf n)+h(\widehat{\mathbf n})
\end{equation}
where
\[
g(\mathbf n)=\frac{1}{\Vert\mathbf n\Vert}\sum_{v\in \mathcal S(\mathfrak p)}\log\vert
1-\phi_{v}(\mathbf n)\vert_v
\]
and
\[
\phi_v(\mathbf n)=\begin{cases}\xi_{\mathfrak p}^{-\mathbf
n}&\mbox{if }\vert\xi_{\mathfrak p}^{\mathbf
n}\vert_v>1;\\
\xi_{\mathfrak p}^{\mathbf n}&\mbox{if
}\vert\xi_{\mathfrak p}^{\mathbf n}\vert_v\leqslant1.\end{cases}
\]
Notice that~$\phi_v(\mathbf n)\neq1$ if~$\mathbf n\neq0$ by the
assumption that the action is mixing.

To establish the lower bound~$C_2>0$ in
Theorem~\ref{theorem:uniformexponentialgrowth}, first note that
equality in~\eqref{equation:generalmoduleperiodicpoints} gives the
expression
\[
\frac{1}{\Vert\mathbf n
\Vert}\log\vert\fix(\alpha_M^{\mathbf n})\vert
\]
as a finite sum of terms of the
form~\eqref{equation:normalizedprimeperiodicpoints} and -- crucially
-- the assumption that~$\dim(X)>0$ means at least one of these terms
arises from a prime~$\mathfrak p$ such that~$\mathbb K(\mathfrak p)$
is an algebraic number field (rather than a function field of
positive characteristic). For the lower bound, it is therefore
enough by~\eqref{equation:generalmoduleperiodicpoints} to consider
only the asymptotic behavior of~$f$ with the assumption that it
arises from such a prime.

We need to know that~$g$ does not make an asymptotic contribution.
It is clear that~$g$ cannot be too large, since
\[
\prod_{v\in \mathcal S(\mathfrak p)}\vert 1-\phi_{v}(\mathbf
n)\vert_v\le2^{\vert\mathcal S(\mathfrak p)\vert}.
\]
On the other hand, it cannot be large and negative for the following
reason. If~$v$ is an infinite (Archimedean) place, then Baker's
theorem~\cite{baker} can be used to find constants~$A,B\geq 0$ such
that
\[
\vert 1-\phi_{v}(\mathbf n)\vert_v\ge\frac{A}{\max\{n_i\}^B}
\]
(see~\cite{MR1461206} or~\cite{emsw} for similar arguments; roughly
speaking, the issue is to bound the proximity of~$\xi_{\mathfrak
p}^{\mathbf n}$ and~$1$ in terms of~$\mathbf n$). If~$v$ is a finite
place, then Yu's~$p$-adic bounds for linear forms in
logarithms~\cite{MR1055245} give a similar lower bound (note
that~$d\ge2$ by assumption). Since~$\mathcal S(\mathfrak p)$ is
finite, this shows that
\begin{equation}\label{equation:applybakerfirst}
g(\mathbf n)\rightarrow0\mbox{ as }\mathbf n\rightarrow\infty.
\end{equation}

Assume, for a contradiction, that~$C_2=0$. Then there is a
sequence~$\mathbf n_j\to\infty$ as~$j\to\infty$ with the property
that
\begin{equation}\label{equation:limitiszero}
\lim_{j\to\infty}f(\mathbf n_j)=0.
\end{equation}
It follows from~\eqref{equation:decomposef}
and~\eqref{equation:applybakerfirst} that
\begin{equation}\label{equation:limitiszeroonsphereforhfunction}
h(\widehat{\mathbf n_j})\rightarrow0\mbox{ as }j\to\infty.
\end{equation}
Now the entropy function~$h$ restricted to the unit sphere~$\mathsf
S_{d-1}$ is a continuous function on a compact set,
so~\eqref{equation:limitiszeroonsphereforhfunction} implies that
there is some~$\mathbf v\in\mathsf S_{d-1}$ for which
\begin{equation}\label{equation:findpointonspherewithzeroentropy}
h(\mathbf v)=0.
\end{equation}
If the ray through~$\mathbf v$ happens to meet a point~$\mathbf
m\in\mathbb Z^d$ then we have an immediate contradiction: the
automorphism~$\alpha_{R_d/\mathfrak p}^{\mathbf m}$ would have zero
entropy, contradicting~\eqref{equation:primemoduleentropy}. In order
to show that any point~$\mathbf v$
with~\eqref{equation:findpointonspherewithzeroentropy} gives a
contradiction we need to make a quantitative version of that
argument.

The explicit formula~\eqref{equation:primemoduleentropy} for~$h$
means there is a list of vectors
\[
\mathbf a_1,\dots,\mathbf a_r\in\mathbb R^d
\]
with the property that for any~$\mathbf u\in\mathbb R^d$,
\[
h(\mathbf u)=\mathbf u\cdot\mathbf a_k,
\]
for some~$k$,~$1\leqslant k\leqslant r$. It
follows that for any~$\mathbf n\in\mathbb Z^d$ and~$\mathbf
u\in\mathbb R^d$,
\[
\vert h(\mathbf n)-h(\mathbf u)\vert\leqslant\max_{1\leqslant
k\leqslant r}\max_{1\leqslant i\leqslant d}\{\vert a_{k,i}\vert\}
\Vert\mathbf n-\mathbf u\Vert,
\]
where~$\mathbf a_k=(a_{k,1},\dots,a_{k,d})$. In particular,
since~$h(\lambda\mathbf v)=0$ for all~$\lambda>0$ we may find a
sequence~$(\mathbf m_j)$ of vectors in~$\mathbb Z^d$
with~$\Vert\mathbf m_j-\lambda_j\mathbf v\Vert\rightarrow0$
as~$j\to\infty$ (for some sequence~$(\lambda_j)$ of scalars) and
hence have
\begin{equation}\label{equation:foracontradiction}
h(\mathbf m_j)\rightarrow0\mbox{ as }j\to\infty.
\end{equation}
On the other hand, the dual group of~$R_d/\mathfrak p$ is connected and
finite-dimensional. Hence, by Yuzvinski{\u\i}'s formula, for
any~$\mathbf n\neq0$,
\[
h({\mathbf n})\geqslant \mah(P)>0
\]
where~$\mah(P)$ denotes the logarithmic Mahler measure of some
polynomial~$P$ of degree no greater than~$\dim_{\mathbb Q}(\mathbb
K(\mathfrak p))$. It follows that there is a constant~$C>0$
(depending only on~$\mathfrak p$) for which
\[
h(\mathbf n)>C>0\mbox{ for any }\mathbf n\neq0
\]
(the existence of a lower bound for the non-zero logarithmic Mahler
measure of polynomials of bounded degree is well-known;
see~\cite{MR0296021} or~\cite{MR1700272} for the background). This
lower bound contradicts~\eqref{equation:foracontradiction},
so~\eqref{equation:limitiszero} is impossible. Thus~$C_2>0$.

The upper bound is clear: the
inequality~\eqref{equation:generalmoduleperiodicpoints} together
with the explicit formula~\eqref{equation:primemoduleperiodicpoints}
gives a uniform constant~$C_1$ with the property that
\[
\vert \fix(\alpha^{\mathbf n})\vert\le C_1^{\max_{1\leqslant
i\leqslant d}\{n_i\}}.
\]
\end{proof}

\begin{remark}
Baker's theorem provides the key estimates in many dynamical
problems; see~\cite{MR1461206} and~\cite{emsw} for examples. As pointed
out by Lind~\cite{MR0684244}, in order to establish the logarithmic
growth rate of periodic points for a quasihyperbolic toral
automorphism (a typical application), sometimes all that is needed
is a weaker and earlier result due to Gel$'$fond~\cite{MR0111736}.
Here we need something closer to the full weight of the theorems of
Baker and Yu, because we are in higher rank.
\end{remark}

\begin{remark}\label{remark:lehmersproblem}
Lehmer's problem~\cite{MR1503118} asks if there is a uniform lower
bound for all positive Mahler measures. As shown by Lind
(see~\cite{MR1062797}) this is equivalent to a uniform lower bound
for the topological entropy of any mixing compact group
automorphism. If `Lehmer's conjecture' that there is such a
bound and that it is attained by the expected polynomial, then the
topological entropy of any compact group automorphism with
positive entropy is at
least
\[
0.162\dots=\mah(x^{10}+x^9-x^7-x^6-x^5-x^4-x^3+x+1).
\]
This does not imply a uniform bound for~$C_2$ in
Theorem~\ref{theorem:uniformexponentialgrowth} because~$C_2$
is influenced by the geometry of the acting group as
well as the collection of maps in its image.
For example, the~$\mathbb Z^2$ action~$\alpha$
corresponding to the
module~$R_2/\langle u_1-2,u_1^ku_2-3\rangle$ has the property
that~$h(\alpha^{(k,1)})=\log3$, so the corresponding
constant~$C_2$ cannot exceed~$\log3/\sqrt{1+k^2}$.
\end{remark}

\end{document}